# Spline surfaces with T-junctions

Kęstutis Karčiauskas  
Vilnius University

Daniele Panozzo  
New York University

Jörg Peters[*]  
University of Florida

May 1, 2017

**Abstract**

*This paper develops a new way to create smooth piecewise polynomial free-form spline surfaces from quad-meshes that include T-junctions, where surface strips start or terminate. All mesh nodes can be interpreted as control points of geometrically-smooth, piecewise polynomials that we call GT-splines. GT-splines are B-spline-like and cover T-junctions by two or four patches of degree bi-4. They complement multi-sided surface constructions in generating free-form surfaces with adaptive layout.*

*Since GT-splines do not require a global coordination of knot intervals, GT-constructions are easy to deploy and can provide smooth surfaces with T-junctions where T-splines can not have a smooth parameterization. GT-constructions display a uniform highlight line distribution on input meshes where alternatives, such as Catmull-Clark subdivision, exhibit oscillations.*

## 1 Introduction

Where strips of surface patches are forced together, it is natural to terminate some; and where strips are stretched out it is natural to spawn additional strips to keep the size and aspect ratio of the patches within bounds. Stopping or initiating surface strips leads to T-junctions where two finer surface pieces meet one coarser piece (see Fig. 1). While in isolation such transitions are easily modeled by smooth hierarchical splines, a full model requires a global coordination of their knot-intervals which is cumbersome or even impossible, as proven by the small mesh in Fig. 2.

The simplest T-junction-configuration, a $\dot{T}$-net (pronounced T1-net), is shown in Fig. 1a: a nominally pentagonal face with exactly one vertex of valence 3 is surrounded by quadrilateral facets. T-junctions allow introducing geometry of higher detail, or to merge two separately-developed spline surfaces. T-junctions also prominently arise when replacing the complex and global constraints of strict quad-meshing [BLP+12, VCD+16] by T-meshes, based on triangle meshes [LRL06, LKH08], curvature directions [ACSD+03, MK04], directional fields [MPKZ10, MPZ14a, PPM+16], optimized for planarity [ZSW10, PW13] or extracted from local parametrizations [RLL+06, JTPSH15].

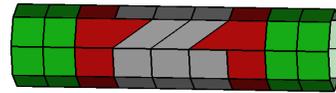

**Figure 2:** *T-splines require that the sum of knot intervals on opposing edges of any face must be equal (Rule 1 of [SZBN03]). This forces the width of the horizontal knot intervals of the grey helical strip to be zero, preventing a smooth T-spline parameterization. (cf. Fig. 13 for a smooth GT-spline surface.)*

**T-junctions and hierarchical splines.** One approach to incorporating T-junctions is hierarchical splines (see e.g. [Kra98, SZBN03, GJS12, DLP13, KXCD15]). Hierarchical splines require that all surface pieces share a single *uv*-parameterization: for any choice of *v*, the *u*-knot intervals must add to the same number; and for any choice of *u*, the *v*-intervals must add to one fixed number. This restriction on the knot sums is natural when refining a single patch. But when the input is a given quadrangulation the knot intervals have to be assigned and hence coordinated. Joining many pieces can then become cumbersome since the local knot intervals have to globally add up to matching sums. Where the mesh is not regular, the local construction in [WZSH11] is therefore only $C^0$ despite generating one order of magnitude more patches than quads. The global construction of [LRL06] introduces additional T-joints and new extraordinary points that are not motivated by geometry but solely by enforcing knot interval constraints (Pre-dating [LZS+12], the construction does not guarantee 'analysis-suitable' knot-distributions.) Akin to [LRL06, Fig.5], the example in Fig. 2 demonstrates that, without modifying the quad-mesh, global coordination is not always possible without loosing smoothness and even continuity of the parameterization. In the example, each red strip forms a bracelet that is half as wide when it comes back to meet up with the starting edge. Since Rule 1 of any T-spline construction according to [SZBN03] mandates that the horizontal knot interval of the red strip be the same where the single-wide edge meets the double-wide edge, the horizontal knot interval of the grey helical strip of patches must be zero. Since, in the example, three consecutive grey horizontal knot intervals are zero, the degree bi-3 spline parameterization is formally $C^{-1}$. That is, the most basic property of T-splines prevents a smooth parameterization for a class of patch layouts that could be hidden in any large scale quad arrangement. The bracelet implies that joining spline surfaces with a *smooth* T-spline parameterization is not always feasible since this, too, requires making equal the knot interval sums. The slightly larger mesh of Fig. 18a demonstrates a more complex incompatibility with any assignment of knot-intervals for smooth T-splines. In summary, while *hierarchical splines are naturally suited for introducing T-junctions in quad meshes*, they are *not naturally suited for generating surfaces from quad meshes with T-junctions*.

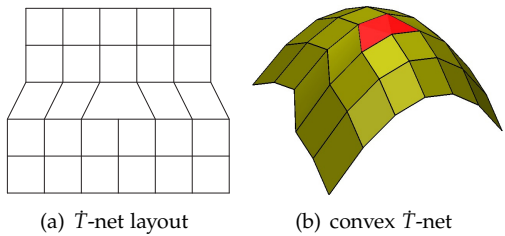

(a) $\dot{T}$-net layout  (b) convex $\dot{T}$-net

**Figure 1:** *A control net with a single isolated T-junction.*

---
[*]Corresponding author, jorg.peters@gmail.com

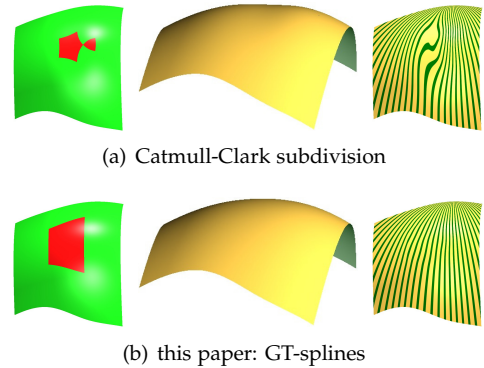

(a) Catmull-Clark subdivision

(b) this paper: GT-splines

**Figure 3:** *Surfaces generated from the convex input $\dot{T}$-net of Fig. 1b. The red regions in (a) represent an infinite sequence of bi-3 patches covering the 5-sided and 3-sided extraordinary Catmull-Clark neighborhoods arising from a $\dot{T}$-net. The red region in (b) consists of two bi-4 patches by which the GT-construction covers the $\dot{T}$-net. Compared to the new construction, (green surface region), Catmull-Clark subdivision produces, in the first two steps, a flattened silhouette and a correspondingly non-uniform highlight line distribution (right).*

**T-junctions and Catmull-Clark subdivision.** A strictly local construction is provided by Catmull-Clark subdivision [CC78]. Here the underlying model is splines with uniform knot spacing and local support. Therefore the parameterizations need not be globally coordinated. However, as Fig. 3a demonstrates for the convex input $\dot{T}$-net of Fig. 1a, the resulting surfaces can be of poor quality: the silhouette dips and rises and the highlight lines oscillate near the T-junction. We note that the oscillations already manifest themselves in the first two subdivision steps and hence rule out Catmull-Clark mesh refinement even as a pre-processor for turning a T-junction into a pair of isolated vertices of valence three and five. (We verified these 'first step artifacts' [ADS11] by replacing, in a sepa-rate computation, the red limit surface in Fig. 3a by a high-quality surface construction.)

**T-junctions and GT-splines.** This paper develops a new local construction, a geometric approach to T-junctions. This GT-construction does not require global knot interval coordination and yields better shape than Catmull-Clark subdivision (see Fig. 3a *vs* 3b). The GT-construction is based on reparameterization, the natural technique for transitioning between unequal parameterizations on opposite sides of a T-junction. Consequently the construction leverages the framework of geometric continuity. For $\dot{T}$-nets, the resulting surfaces

- consist of a frame of bi-3 (bi-cubic) patches filled by two patches of degree bi-4 (this 'cap' is red in Fig. 1b).
- The bi-4 cap is internally smooth and joins the bi-cubic patches with tangent continuity ($G^1$)
- The bi-4 cap yields good highlight line distributions on all of a large number of challenging input meshes.
- The Appendix provides simple explicit formulas for all relevant Bernstein Bézier coefficients of the GT-construction in terms of the local $\dot{T}$-net.
- GT-splines complement few-piece polynomial constructions such as [KNP16] covering extraordinary points to model smooth free-form surfaces of maximal degree bi-4.

We also present two variants of the GT-construction, to generate smooth caps for multiple T-junctions within one facet whose extensions cross or are parallel.

**Overview.** Section 2 reviews basic concepts of the construction of smooth surfaces. Section 3 defines a bi-3 frame of patches that transitions to the surrounding surface. Section 4 describes the $G^1$ bi-4 GT-construction of the cap for $\dot{T}$-nets. Section 5 and Section 6 develop caps for cases of adjacent T-junctions: for two T-junctions opposite one and for two T-junctions with crossing directions. Section 7 compares the constructions for challenging input data, explains the choices taken along the way and lists the limitations. Section 8 shows how GT-splines collaborate with algorithms in the literature to smoothly cover multi-sided neighborhoods by an overall smooth bi-4 surface.

## 2 Definitions and Setup

The GT-splines are a collection of tensor-product patches in Bernstein-Bézier form (BB-form; see e.g. [Far88]):

$$\mathbf{f}(u,v) := \sum_{i=0}^{d} \sum_{j=0}^{d} \mathbf{f}_{ij} B_i^d(u) B_j^d(v), \quad (u,v) \in [0..1]^2,$$

where $B_k^d(t) := \binom{d}{k}(1-t)^{d-k}t^k$ are the BB polynomials of degree $d$ and $\mathbf{f}_{ij}$ are the BB-coefficients. Adjacent patches join with $G^k$ continuity if their $k$th-order jets (one-sided Taylor expansion) match along their common boundary after a change of variables $\rho$. This characterization is equivalent to formulations of $C^k$ continuity of manifolds in terms of charts, see e.g. [Pet02]. We use the succinct characterization that *two surface pieces $\tilde{\mathbf{f}}$ and $\mathbf{f}$ sharing a boundary curve $\mathbf{e}$ join $G^1$ if there is a suitably oriented and non-singular reparameterization $\rho : \mathbb{R}^2 \to \mathbb{R}^2$ so that the jets $\partial^k \tilde{\mathbf{f}}$ and $\partial^k (\mathbf{f} \circ \rho)$ agree along $\mathbf{e}$ for $k = 0, 1$.* Although $\rho$ is just a change of variables, its choice is crucial for the properties of the resulting surface. Throughout, we will choose $\mathbf{e}$ to correspond to the patch parameters $(u, 0 = v)$. Then the relevant Taylor expansion of the reparameterization $\rho$ with respect to $v$ is $\rho := (u + b(u)v, a(u)v)$ and the chain rule of differentiation yields the $G^1$ constraints

$$\partial_v \tilde{\mathbf{f}} - a\partial_v \mathbf{f} - b\partial_u \mathbf{f} = 0. \tag{1}$$

If $\tilde{\mathbf{f}}$ and $\mathbf{f}$ are polynomials then where $a$ and $b$ are rational functions of degree bounded in terms of the degree of $\tilde{\mathbf{f}}$ and $\mathbf{f}$ [Pet91].

$C^1$ continuity of the splines over non-uniform knot sequences can be recast as $G^1$ continuity of patches defined over unit domains (see e.g. [KP11]). For example, if $\mathbf{f}$ and $\tilde{\mathbf{f}}$ are consecutive curve segments originally associated with intervals $[-1, 0]$ and $[0, \frac{1}{2}]$ of a $C^1$ spline with knot sequence $\{\ldots, -1, 0, \frac{1}{2}, \ldots\}$ then both $\mathbf{f}$ and $\tilde{\mathbf{f}}$ can be newly defined, each on the interval $[0, 1]$,

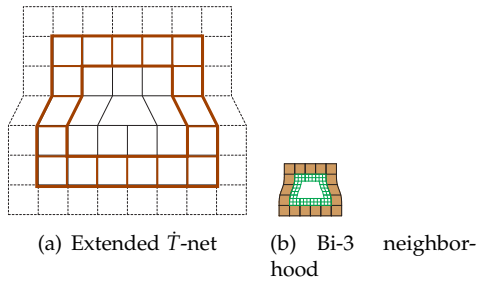

(a) Extended $\dot{T}$-net  (b) Bi-3 neighborhood

**Figure 4:** *An isolated T-junction in an extended $\dot{T}$-net (a) provides (b) a bi-3 neighborhood (solid, for context only) and a $C^2$-prolongation in BB-form (inner green mesh) that is the only part used for the GT-construction.*

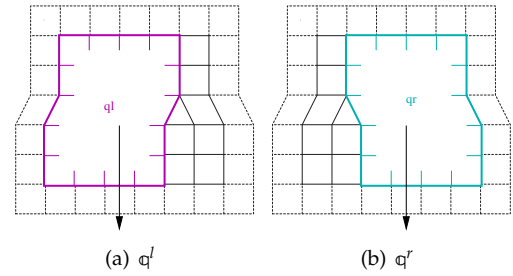

(a) $\mathfrak{q}^l$  (b) $\mathfrak{q}^r$

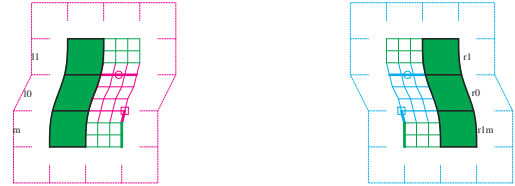

(c) Interpretation as BB-patches of degree bi-3, overlaid in (e)

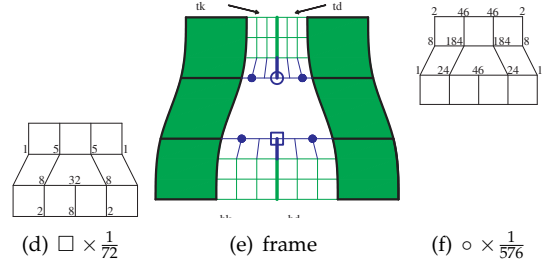

(d) $\square \times \frac{1}{72}$  (e) frame  (f) $\circ \times \frac{1}{576}$

**Figure 5:** *Constructing the frame. The regular left (a) and right (b) control nets obtained by re-connecting the nodes of the $\dot{T}$-net define six bi-3 patches each. (e) Completion of the frame. (d,f) Stencils of the points $\mathbf{q}^{b,-1}_\square = \mathbf{q}^{b,1}_\square$ marked $\square$ and $\mathbf{q}^{t,-1}_\circ = \mathbf{q}^{t,1}_\circ$ marked by $\circ$. in (e).*

and they then join as $\beta \partial_u \mathbf{f}(1) = \partial_u \tilde{\mathbf{f}}(0)$ with $\beta := 1/2$. When we want to point out that surface patches are, in one variable, related by the identity and, in the other, by $C^1$ continuity over non-uniform knot sequence, we refer to the transition as: $C^1$ *with parameter $\beta$*.

Our main construction focusses on $\dot{T}$-nets that consist, as shown in Fig. 1, of quadrilaterals and one nominally five-sided facet. For context and exposition, we can extend the $\dot{T}$-net by one layer of quadrilaterals (see Fig. 4a). This allows applying, away from the five-sided facet, the well-known bi-cubic (bi-3) B-spline to BB-form conversion rules (see e.g. [Far88]). The resulting $C^2$ bi-3 spline neighborhood is colored brown in Fig. 4b. The smaller $\dot{T}$-net ( Fig. 1 ) provides a *tensor-border of degree* 3 *and depth* 2, the $C^2$-prolongation of second-order Hermite data shown as a green net of BB-coefficients in Fig. 4b. Given the tensor-border, the $\dot{T}$-net interior of $4 + 4 + 5 + 5$ control points (see the stencils in Fig. 5d,f) provides all the information for the GT-construction. All formulas (stencils) of the bi-4 cap in terms of the $\dot{T}$-net interior are provided in the Appendix.

# 3 Construction of a bi-3 frame of patches for capping a $\dot{T}$-net

The surface corresponding to a $\dot{T}$-net will consist of a frame of bi-3 patches and a central cap consisting of two bi-4 patches. This section builds the frame. For $i = -1, 0, 1$ the frame has left patches $\mathbf{q}^{l,i}$, right patches $\mathbf{q}^{r,i}$, (Fig. 5c) and for $j \in \{-1, 1\}$ top patches $\mathbf{q}^{t,j}$ and bottom patches $\mathbf{q}^{b,j}$ (Fig. 5e). The ribbon is derived by *re-connecting* the nodes of the $\dot{T}$-net to form two regular nets, $\mathfrak{q}^l$ from the left (see Fig. 5a) and $\mathfrak{q}^r$ from the right (see Fig. 5b).

We now interpret $\mathbb{q}^l$ and $\mathbb{q}^r$ as bi-3 B-spline control nets and convert them to BB-form. As illustrated in Fig. 5c the so-derived patches agree – except for the lower boundary curves of the top patches and the upper boundary of the bottom patches. The bottom curves overlap only at their endpoints, marked as a magenta and a cyan box in Fig. 5c. A new common BB-control point, also marked as a blue box in Fig. 5e, is chosen to be the average of these two candidates: Fig. 5d displays the explicit stencil for the common point $\mathbf{q}_\square := (\mathbf{q}_\square^{b,-1} + \mathbf{q}_\square^{b,1})/2$ obtained by this averaging. The two direct neighbors of $\mathbf{q}_\square$ are chosen to make the combined boundary curve $C^2$. (The blue disks are defined by the $C^1$-prolongation of $\mathbf{q}^{l,-1}$ and $\mathbf{q}^{r,-1}$). The top patches are subdivided at their midpoint and the resulting overlap is then treated like that of the bottom: Fig. 5f displays the explicit stencil for the common point $\mathbf{q}_\circ$ obtained by this averaging. Although this split into $\mathbf{q}^{t,-1}$ and $\mathbf{q}^{t,1}$ serves only to accommodate the lower boundary of the top patches, extensive experiments show this split to be critical for achieving good shape (see e.g. Fig. 12). The resulting frame of bi-3 patches is $C^2$ except along the four *hv-curves*, the red curves in Fig. 6c between the horizontal and the vertical strips of the frame: $\mathbf{q}^{l,-1}$ to $\mathbf{q}^{b,-1}$, $\mathbf{q}^{l,1}$ to $\mathbf{q}^{t,-1}$, $\mathbf{q}^{r,-1}$ to $\mathbf{q}^{b,1}$, $\mathbf{q}^{r,1}$ to $\mathbf{q}^{t,1}$. Across the hv-curves the continuity is $C^1$. Since the construction did not change the BB-coefficients derived from $\mathbb{q}^l$ and $\mathbb{q}^r$ that match the tensor-border (Fig. 4b), the frame joins $C^2$ with the splines surrounding it.

## 4 Central $\dot{T}$-net cap construction

To complete the surface, we construct a central cap (red in Fig. 6) of degree bi-4 that fills the frame so that all transitions are at least $G^1$. The bi-4 cap consists of two patches $\mathbf{p}^l$, $\mathbf{p}^r$; see Fig. 6c. Since the construction of $\mathbf{p}^r$ mirrors that of $\mathbf{p}^l$, we discuss only $\mathbf{p}^l$. Fig. 6a shows the $C^1$-prolongations t of $\mathbf{q}^{t,-1}$ and b of $\mathbf{q}^{b,-1}$ in black and l of $\mathbf{q}^{l,0}$ in green. While b is consistent with l, the prolongations t and l are inconsistent (due to the split of the top patch when constructing the frame). Since we reparameterize l linearly (to minimize the final patch degree) to match t, we also need to reparameterize b after all. Together, the choice of parameterizations in Eq. (1) are

$$
\begin{array}{rcc}
& a(u) := & b(u) := \\
\text{left(l)}: & 1 - \frac{u}{2} & 0 \\
\text{top(t)}: & 1 & (1-u)u \\
\text{bottom(b)}: & 1 & -\frac{1}{2}(1-u)u.
\end{array} \quad (2)
$$

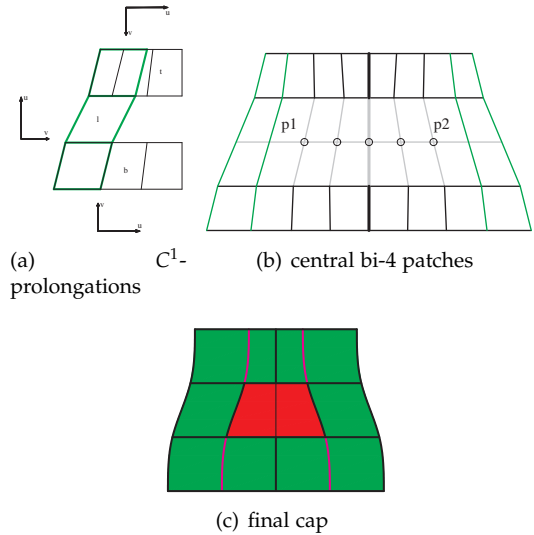

(a) $C^1$-prolongations   (b) central bi-4 patches

(c) final cap

**Figure 6:** *Construction of the bi-4 cap. (a) the mismatch of the $C^1$-prolongations is resolved by reparameterizing them. (b) Interior coefficients minimize the distance to bi-cubics. (c) final layout: the bi-3 patches are $C^1$-connected across the (red) hv-curves.*

The interior BB-coefficients (circles in Fig. 6b) of the bi-4 cap are determined so that columns of BB-coefficients form degree-raised curves of true degree 3. By construction, see Fig. 6c, the red bi-4 cap is internally $C^1$ and joins with $G^1$-continuity to the green frame.

# 5 Caps for parallel T-junctions

Configurations with multiple T-junctions can in principle be locally re-meshed to separate them into isolated $\dot{T}$-nets. For completeness, and to compare what surface quality can be achieved, we investigate configurations where two T-junctions face another, as shown in Fig. 7b. We call the configuration a $\ddot{T}$-net (pronounced T3-net). It has one nominally 7-sided face. Such T-junctions can arise, for example, from configuration Fig. 7a by removing the two red edges of a triangle attached to a point of valence 5. (Asymmetrically removing one yields a mesh with one T-junction as in Fig. 7c). Applying Catmull-Clark subdivision to $\ddot{T}$-junctions leads to poor surfaces. As for $\dot{T}$-nets, first a bi-3 frame is constructed. Fig. 7e,f provide the stencils for the points marked □ and ○. Except across the red hv-curves between the horizontal and the vertical strips, the frame is $C^2$. In the spirit of the bi-4 $\dot{T}$-net construction, a bi-4 cap for the $\ddot{T}$-net is constructed by subdividing, in the ratio shown in Fig. 8a, the $C^1$-prolongation of the frame from the top; and then evenly splitting the middle prolongation of the bottom. Both the top and bottom prolongations are $C^1$-connected (in the horizontal direction) with the same continuity parameters from left to right: $\beta = \frac{1}{2}$, $\beta = 1$, $\beta = 2$. This implies a continuity pa-rameter of $\beta = \frac{2}{3}$ across the top hv-curves. For the prolongation 1 of $\mathbf{q}^{l,0}$ (see Fig. 8b) to match the split prolongation of the top, the data are reparameterized according to

$$\begin{array}{ccc} & a(u) := & b(u) := \\ \text{left}(1): & 1 - u/3 & 0 \\ \text{top}(\text{left}): & 1 & \frac{1}{2}(1-u)u \\ \text{bottom}(\text{left}): & 1 & -\frac{1}{3}(1-u)u. \end{array} \quad (3)$$

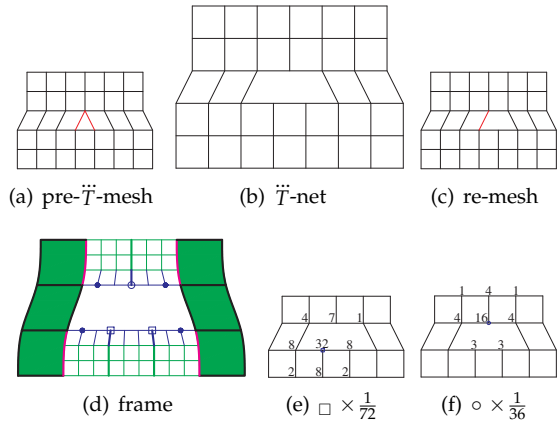

(a) pre-$\ddot{T}$-mesh  (b) $\ddot{T}$-net  (c) re-mesh

(d) frame  (e) □ × $\frac{1}{72}$  (f) ○ × $\frac{1}{36}$

**Figure 7:** *Construction of the frame for a $\ddot{T}$-junction.*

(For the right top and bottom reparameterizations, $b(u)$ is negated.) The remaining (circled) BB-coefficients in Fig. 8 are chosen to make their columns have actual degree 3.

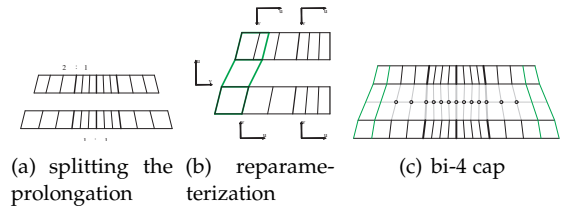

(a) splitting the prolongation  (b) reparameterization  (c) bi-4 cap

**Figure 8:** *bi-4 cap for a $\ddot{T}$-net.*

# 6 Caps for two T-junctions in crossing directions

We also investigate configurations with two T-junctions as shown in Fig. 9a. We call the configuration a $\ddot{T}$-net (pronounced T2-net). It has one nominally 6-sided face. Note that such configurations are explicitly excluded in dyadic T-meshes [KBZ15].

Capping $\ddot{T}$-nets with good surface quality is similar but more challenging than the earlier constructions. As before, B-spline to Bézier conversion yields the green Bézier control points in Fig. 9b, now with both the top and the right side patches subdivided.

The corner points (blue boxes in Fig. 9b) and the middle, blue circle coefficients of innermost (blue) boundary curves of the frame are determined only in the last step of the construction. The direct neighbors of the corner points are chosen so that adjacent bi-3 patches of the frame connect $C^1$ (with ratios 1,1 lower left, 1/2,1/2

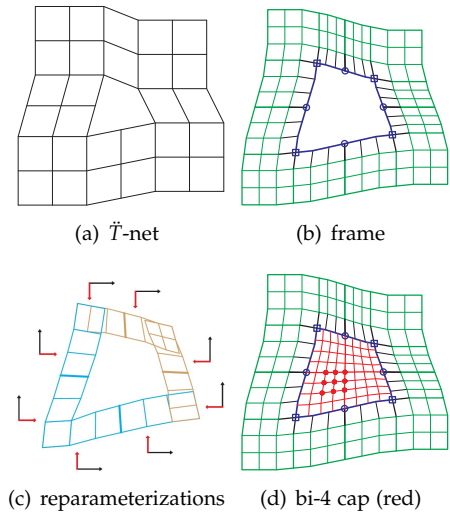

(a) $\ddot{T}$-net  (b) frame

(c) reparameterizations  (d) bi-4 cap (red)

**Figure 9:** *Bi-4 cap for a $\ddot{T}$-configuration. The red axes in (c) indicate the v-parameter.*

upper right and 1/2,1 otherwise); the direct neighbors of the middle points are determined so that each boundary curve is internally $C^2$.

The $C^1$-prolongations are not compatible at the corners (see Fig. 9c). To make them compatible they are re-parameterized with

$$
\begin{array}{lll}
 & a(u) := & b(u) := \\
\text{top} - \text{left}: & 1 - \frac{u}{4} & \frac{1}{2}(1-u)u \\
\text{top} - \text{right}: & \frac{3}{4} - \frac{u}{4} & -\frac{1}{2}(1-u)u \\
\text{bottom} - \text{left}: & 1 - \frac{u}{4} & -\frac{1}{4}(1-u)u \\
\text{bottom} - \text{right}: & \frac{3}{4} - \frac{u}{4} & \frac{1}{4}(1-u)u.
\end{array}
$$

This list of the reparameterizations is complete due to the (combinatorially) diagonal symmetry (see the local coordinate systems in Fig. 9c). The reparameterized tensor-border is of degree 4 and depth 1 and ensures $G^1$ continuity of the central cap with the frame. Choosing to join $C^2$ the four $3 \times 3$ groups of interior BB-coefficients leaves free one group shown as red disks. Finally, we minimize, over all 11 bi-3 patches of the frame and the 4 bi-4 patches of the central cap, the functional $\mathcal{F}_3$ where

$$\mathcal{F}_\kappa f := \int_0^1 \int_0^1 \sum_{i+j=\kappa, i,j\geq 0} \frac{\kappa!}{i!j!}(\partial_s^i \partial_t^j f(s,t))^2 dsdt.$$

We minimize the sum with respect to 17 unknown coefficients: 4 corner (blue box), 4 mid-edge (blue circle) and 9 inner ones (red disks in Fig. 9d). In the implementation, these 17 coefficients enter as affine combinations of $\ddot{T}$-net points with pre-computed coefficients.

The choice $\mathcal{F}_3$ is the result of testing a series of input meshes including the challenging elliptic configuration in Fig. 10. Fig. 11 confirms that this choice also works well for a wave-like input mesh.

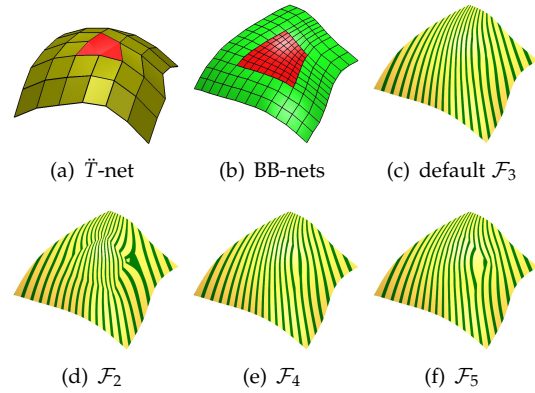

(a) $\ddot{T}$-net  (b) BB-nets  (c) default $\mathcal{F}_3$

(d) $\mathcal{F}_2$  (e) $\mathcal{F}_4$  (f) $\mathcal{F}_5$

**Figure 10:** *The surfaces obtained by minimizing functionals. (a) convex $\ddot{T}$-net. (b) BB-coefficients of the frame (green) and the central bi-4 cap (red). (c-f) Highlight lines.*

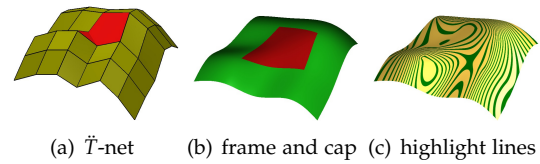

(a) $\ddot{T}$-net  (b) frame and cap  (c) highlight lines

**Figure 11:** *A $\ddot{T}$-net ($\mathcal{F}_3$) surface for a non-convex input mesh.*

## 7 Discussion and Comparison

### 7.1 The frame construction: to split or not to split

Note the flaws in the highlight line distribution of the surface in Fig. 12b. For the same input mesh Fig. 1b, the highlight line distribution of Fig. 12b is much worse than that of the GT-construction Fig. 3b. Despite being formally smooth, surfaces generated by not subdividing the frame patch opposite the T-junction have poor highlight lines. Evidently splitting, though not required by formal smoothness contraints, improves shape quality.

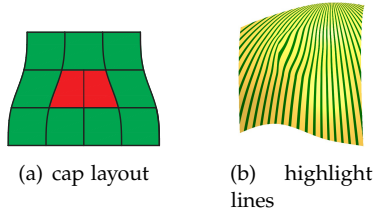

(a) cap layout     (b) highlight lines

**Figure 12:** *An alternative construction that does not split the patch opposite to the T-junction leads to a poor highlight line distribution.*

## 7.2 Comparing GT-splines to T-splines

In general a comparison to hierarchical splines does not make sense since we proved at the outset that not all meshes with T-junctions admit smooth T-splines. For one such configuration, the bracelet mesh of Fig. 2 (that does not admit a $C^1$ T-spline), Fig. 13 demonstrates that applying GT-splines yields a bi-4 surface with an excellent highlight line distribution.

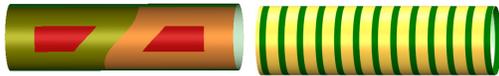

**Figure 13:** *The GT-spline surface for the T-mesh of Fig. 2 that does not admit a $C^1$ T-spline. The bi-4 cap in red; the right image shows highlight lines.*

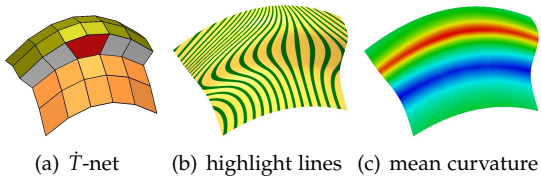

(a) $\dot{T}$-net     (b) highlight lines     (c) mean curvature

**Figure 14:** *Two regular meshes are merged with a T-junction.*

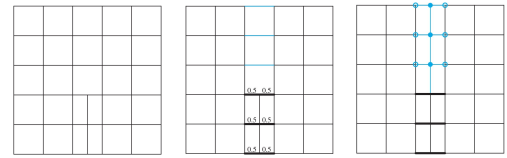

(a) GT-spline input     (b) T-spline input     (c) B-spline input

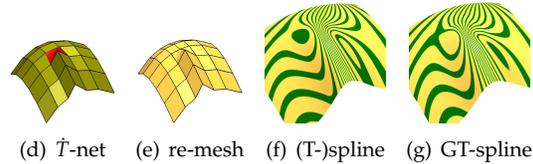

(d) $\dot{T}$-net     (e) re-mesh     (f) (T-)spline     (g) GT-spline

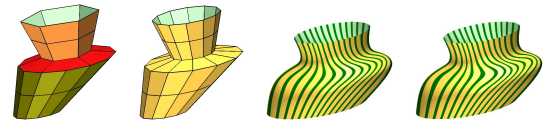

(h) input mesh     (i) re-mesh     (j) (T-) spline     (k) GT-spline

**Figure 15:** *Comparisons to T-splines and tensor-product splines. Since, in these examples, the T- and B-spline surfaces coincide, we refer to them as (T-)spline surfaces.*

Fig. 14 shows the case of two regular meshes of different quad-patch count joined via a T-junction. T-splines require forming a common parameter domain, that, while easy for simple meshes, is impossible for more complex meshes such as Fig. 18. The result of directly applying a GT-spline is displayed in Fig. 14b,c. To establish an upper bound on the quality of the highlight line distribution, Fig. 15 compares GT-splines to T-splines in the form of standard bi-cubic tensor-product splines. The top row of Fig. 15 shows in order (a) the $\dot{T}$-net input to the GT-construction, (b) the geometrically identical T-mesh, with all knot intervals are 1 except for 0.5 on the thick edges; and (c) the mesh resulting from splitting the cyan knot intervals yielding new (cyan disks) and moved (cyan circles) control points. Since the T-spline surface of input mesh (b) equals the non-uniform $C^2$ tensor-product B-spline surface of input mesh (c), it suffices to discuss the case (c) in the following. For a concrete comparison in this setting, the mesh in Fig. 15d is the geometric input both for the GT-spline and for the T-spline, whereas Fig. 15e is the tensor-product B-spline input mesh. As explained above, the knot-intervals are chosen so that the B-spline surface equals the T-spline surface. Fig. 15f,g show the highlight line distribution on the surfaces.

In Fig. 15h two regular meshes are connected by a whole ring of T-junctions. This is the (geometric) input mesh both for the GT-spline and for the T-spline, but the knot-intervals differ: for the GT-spline they can all be chosen equal, while the horizontal T-spline knot intervals of bottom mesh are 0.5 when those of the top mesh are 1. Fig. 15i shows the tensor-product mesh yielding the $C^2$ surface Fig. 15j that coincides with T-spline surface from mesh Fig. 15h. Fig. 15k shows the GT-spline surface obtained directly from Fig. 15h. The highlight lines of the GT-construction are very similar to those of the (T-)spline construction even though the GT-spline is placed at a disadvantage by capping T-junctions while the (T-)spline can take advantage of a tensor-product mesh.

Given the similarity in shape, it is important to recall the essential difference between T-splines and GT-splines. The simple knot intervals used *locally* above can lead to invalid knot-intervals when considering larger meshes: T-splines require a global coordination of knot sequences. Such global coordination is not always possible or may require complex re-meshing. The GT-construction is local, sidestepping the need for global coordination, and producing surfaces of comparable quality.

## 7.3 Separation and Re-meshing

The GT-constructions in Section 4, 5 and 6 follow a common pattern, of adjusting the bi-3 frame and then forming a central cap of degree bi-4. The minimal submesh required for the GT-constructions is called 'net'. Only the outer boundaries of the frame may have *irregular nodes*, where $n \neq 4$ quads meet, or T-junctions. Guaranteeing such separation may not be easy in general (see Limitations below) but often local adjustments can be made. While Section 6 demonstrated that $\ddot{T}$-nets can yield bi-4 surfaces of good quality, better surfaces are often obtained by locally re-connecting the mesh points to isolate the T-junctions as in Fig. 16c. Reconnecting with an even larger footprint for a more symmetric reconnection as in Fig. 16d can further improve the highlight line distribution, at the cost of higher construction complexity.

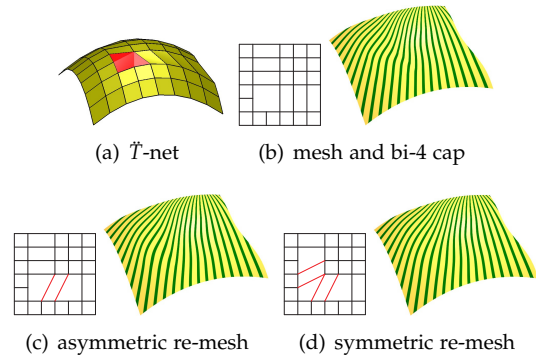

(a) $\ddot{T}$-net  (b) mesh and bi-4 cap

(c) asymmetric re-mesh  (d) symmetric re-mesh

**Figure 16:** *Effect of re-meshing on the resulting bi-4 cap.*

## 7.4 Limitations

Quadrangulations can contain closely packed T-junctions and irregular points. As our initial example demonstrated, Catmull-Clark refinement is not a good way to separate T-junctions: not only does subdivision increase the number of patches, but, more importantly, it can negatively impact the surface quality. T-mesh subdivision [KBZ15] can also start with irregular points adjacent to T-junctions but does not separate T-junctions from irregular points. Therefore, it too cannot be used for pre-processing.

Although re-meshing can reduce many configurations to the three standard T-nets, GT-constructions are not expected to work with arbitrary T-junction distributions. Here, the quad-meshing algorithm or the designer have to enforce some discipline, already to obtain good shape. Many-sided facets with T-junctions as generated by [ACSD+03] or motorcycle graphs [EGKT08, MPZ14b] are outside the scope of GT-constructions.

In many cases, our approach can allow for tighter packing of T-junctions (for example as in Fig. 15,*bottom* row and Fig. 18) and irregularities. However, this paper does not attempt to provide a set of recipes for arbitrarily complex T-junctions. On one hand, local re-connection can often reduce the situation to a collection of $\dot{T}$-net, $\ddot{T}$-nets or $\dddot{T}$-nets, but a principled prescription for such quad-re-meshing is outside the scope. On the other hand, our experiments with complex configurations show that keeping the complexity of capping T-junctions to a minimum results in better shape.

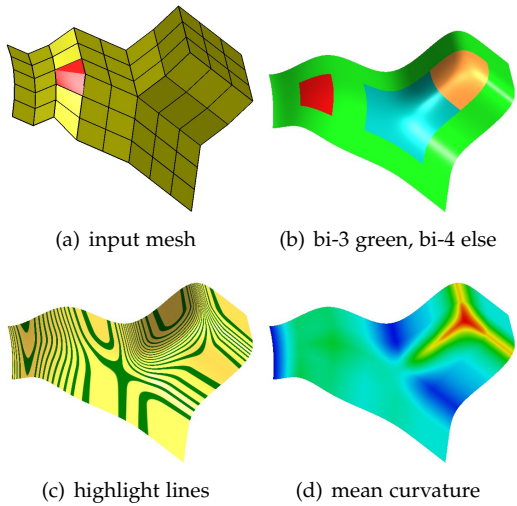

(a) input mesh  (b) bi-3 green, bi-4 else

(c) highlight lines  (d) mean curvature

**Figure 17:** *Mesh and bi-4 surface combining a T-junction with irregular regions of valence* 3 *and* 5 *where [KNP16] is applied.*

# 8 Collaboration with multi-sided caps

Fig. 17 demonstrates that sufficiently isolated caps for T-junctions co-exist without problems with irregular vertices where the surface caps are also of degree bi-4 when we apply [KNP16]. Again, tighter configurations are possible, but may reduce surface quality while increasing the complexity of implementation.

Fig. 18 presents another free-form design that challenges algorithms that require globally consistent knot intervals. As in Fig. 2, enforcing

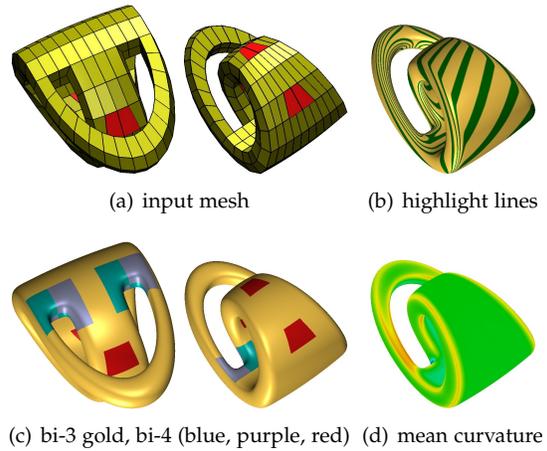

(a) input mesh  (b) highlight lines

(c) bi-3 gold, bi-4 (blue, purple, red)  (d) mean curvature

**Figure 18:** *Mesh and $G^1$ surface including horizontally paired $\dot{T}$-nets (red) and irregular neighborhoods (blue,purple) treated with [KNP16]. This mesh does not admit a globally consistent (non-zero) knot interval assignment for smooth T-splines.*

T-spline Rule 1 yields zero knot intervals, now also at extraordinary points.

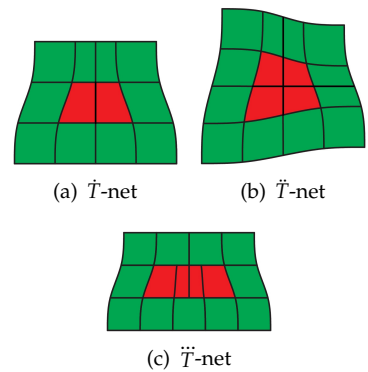

(a) $\dot{T}$-net  (b) $\ddot{T}$-net

(c) $\dddot{T}$-net

**Figure 19:** *Distribution of the bi-degree 4 patches (red) for the three basic configurations. The frames (green) are of degree bi-3.*

# 9 Conclusion

The paper introduced a construction of surface caps for merging and spreading feature lines via T-junctions. In the default $\dot{T}$-case, the caps consist of two surface patches of degree bi-4

Fig. 19, otherwise of four bi-4 patches. The surrounding bi-3 patches are only perturbed where they join the central cap. Since the approach is based on geometric continuity, it does not require non-local coordination of knot intervals and is not restricted to graphs of functions but applies to general manifolds.

# Appendix

Due to symmetry and $C^2$ continuity of the bi-3 frame with the surrounding surface, the $\dot{T}$-net-construction is completely defined by the following formulas for one (here the 'left') bi-4 patch. As illustrated for two of the coefficients in Fig. 5d,f, each of the $5 \times 5$ BB-coefficient is a linear combination of the inner $\dot{T}$-net nodes with two rows of four weights and two rows of five weights. In Table I, the stencil for each BB-coefficient, is scaled $\times\ 144$ and placed into brackets.

$$\begin{bmatrix} 4 & 16 & 4 & 0 & 0 \\ 16 & 64 & 16 & 0 & 0 \\ 4 & 16 & 4 & 0 & 0 \\ 0 & 0 & 0 & 0 & 0 \\ 0 & 0 & 0 & 0 & 0 \end{bmatrix} \begin{bmatrix} \frac{5}{2} & 16 & \frac{11}{2} & 0 & 0 \\ 10 & 64 & 22 & 0 & 0 \\ 0 & 0 & 0 & 0 & 0 \\ 0 & 0 & 0 & 0 & 0 \\ 0 & 0 & 0 & 0 & 0 \end{bmatrix} \begin{bmatrix} 0 & \frac{3}{8} & 15 & \frac{15}{2} & 0 \\ 0 & 6 & 60 & 30 & 0 \\ 0 & 13 & \frac{35}{4} & 0 & \frac{-1}{8} \\ 0 & 0 & 0 & 0 & 0 \\ 0 & 0 & 0 & 0 & 0 \end{bmatrix} \begin{bmatrix} \frac{7}{8} & \frac{107}{2} & \frac{77}{2} & \frac{1}{8} & \frac{-1}{8} \\ \frac{5}{2} & 9 & 0 & 3 & 0 \\ \frac{23}{8} & \frac{413}{8} & \frac{297}{8} & \frac{231}{8} & \frac{-1}{8} \\ \frac{161}{8} & \frac{323}{8} & 6 & 0 & 0 \\ 0 & 0 & 0 & 0 & 0 \end{bmatrix} \begin{bmatrix} \frac{1}{2} & \frac{23}{8} & \frac{46}{2} & \frac{23}{8} & \frac{1}{4} \\ 2 & 6 & 6 & 6 & 2 \\ \frac{21}{8} & \frac{46}{2} & \frac{161}{8} & \frac{46}{2} & \frac{21}{8} \\ 0 & \frac{23}{8} & 14 & \frac{23}{8} & 0 \\ 0 & \frac{7}{16} & 0 & \frac{7}{16} & 0 \end{bmatrix}$$

$$\begin{bmatrix} 1 & 4 & 1 & 0 & 0 \\ 16 & 64 & 16 & 0 & 0 \\ 28 & 7 & 0 & 0 & 0 \\ 0 & 0 & 0 & 0 & 0 \\ 0 & 0 & 0 & 0 & 0 \end{bmatrix} \begin{bmatrix} \frac{1}{4} & 4 & \frac{7}{2} & 0 & 0 \\ \frac{7}{2} & 64 & \frac{47}{2} & 0 & 0 \\ 28 & 10 & 0 & 0 & 0 \\ 0 & 0 & 0 & 0 & 0 \\ 0 & 0 & 0 & 0 & 0 \end{bmatrix} \begin{bmatrix} \frac{1}{8} & 5 & \frac{13}{4} & 0 & 0 \\ 0 & \frac{85}{4} & 58 & \frac{21}{8} & 0 \\ \frac{67}{32} & 116 & 33 & \frac{9}{4} & \frac{-9}{32} \\ 0 & 43 & 0 & 0 & 0 \\ 0 & 0 & 0 & 0 & 0 \end{bmatrix} \begin{bmatrix} \frac{1}{8} & \frac{23}{8} & \frac{23}{8} & \frac{23}{8} & \frac{-1}{8} \\ \frac{25}{8} & 0 & \frac{-5}{16} & 0 & 0 \\ 0 & \frac{47}{2} & \frac{297}{8} & \frac{5}{16} & \frac{1}{16} \\ 15 & 0 & 0 & \frac{5}{2} & 0 \\ 0 & 0 & 0 & \frac{-1}{16} & 0 \end{bmatrix} \begin{bmatrix} \frac{1}{2} & 2 & \frac{23}{8} & 6 & \frac{1}{4} \\ 3 & 0 & 33 & 0 & 3 \\ 14 & \frac{87}{2} & 33 & 14 & \frac{1}{4} \\ 0 & \frac{35}{2} & 0 & \frac{35}{2} & 0 \\ 0 & 0 & 3 & 0 & 0 \end{bmatrix}$$

$$\begin{bmatrix} 12 & 0 & 0 & 0 & 0 \\ 0 & 48 & 12 & 0 & 0 \\ 48 & 12 & 0 & 0 & 0 \\ 0 & 0 & 0 & 0 & 0 \\ 0 & 0 & 0 & 0 & 0 \end{bmatrix} \begin{bmatrix} 0 & \frac{9}{2} & 0 & 0 & 0 \\ 6 & 48 & \frac{39}{2} & 0 & 0 \\ 0 & 18 & 0 & 0 & 0 \\ 0 & 0 & 0 & 0 & 0 \\ 0 & 0 & 0 & 0 & 0 \end{bmatrix} \begin{bmatrix} \frac{7}{6} & \frac{-1}{6} & \frac{-1}{3} & \frac{1}{2} & 0 \\ \frac{19}{6} & 43 & 25 & \frac{5}{6} & 0 \\ \frac{116}{3} & 31 & \frac{4}{3} & 0 & \frac{-1}{6} \\ 0 & \frac{-3}{2} & 0 & 0 & 0 \\ 0 & 0 & 0 & 0 & 0 \end{bmatrix} \begin{bmatrix} \frac{9}{16} & 0 & 0 & 0 & 0 \\ 0 & \frac{7}{2} & \frac{325}{16} & 0 & \frac{7}{2} \\ 0 & \frac{59}{2} & \frac{235}{16} & \frac{5}{2} & 0 \\ 0 & \frac{17}{8} & 4 & \frac{-1}{8} & 0 \\ 0 & 0 & 0 & 0 & 0 \end{bmatrix} \begin{bmatrix} 0 & 3 & 0 & \frac{21}{2} & 0 \\ 14 & \frac{33}{2} & 33 & 3 & 0 \\ 0 & \frac{87}{2} & 14 & \frac{1}{4} & 0 \\ 0 & \frac{35}{2} & 0 & \frac{35}{2} & 0 \\ 0 & 4 & 0 & 4 & 0 \end{bmatrix}$$

$$\begin{bmatrix} 0 & 0 & 0 & 0 & 0 \\ 16 & 64 & 28 & 0 & 0 \\ 1 & 4 & 1 & 0 & 0 \\ 0 & 0 & 0 & 0 & 0 \\ 0 & 0 & 0 & 0 & 0 \end{bmatrix} \begin{bmatrix} 0 & \frac{17}{2} & 28 & 0 & 0 \\ \frac{11}{2} & 64 & \frac{53}{2} & 0 & 0 \\ \frac{5}{8} & 4 & \frac{11}{8} & 0 & 0 \\ 0 & 0 & 0 & 0 & 0 \\ 0 & 0 & 0 & 0 & 0 \end{bmatrix} \begin{bmatrix} 0 & 2 & 25 & 14 & \frac{1}{2} \\ 0 & 4 & 52 & 44 & 2 \\ 0 & 4 & 0 & 0 & 0 \\ 0 & 0 & 0 & 0 & 0 \\ 0 & 0 & 0 & 0 & 0 \end{bmatrix} \begin{bmatrix} 0 & \frac{7}{2} & 28 & 0 & 0 \\ 0 & \frac{59}{2} & \frac{17}{8} & 64 & 4 \\ 0 & \frac{17}{8} & 0 & 4 & 0 \\ 0 & 0 & 0 & 0 & 0 \\ 0 & 0 & 0 & 0 & 0 \end{bmatrix} \begin{bmatrix} 0 & 0 & 0 & 0 & 0 \\ \frac{7}{2} & \frac{35}{2} & \frac{35}{2} & \frac{7}{2} & 0 \\ 16 & \frac{87}{2} & 33 & 16 & 1 \\ 4 & \frac{35}{2} & 3 & 4 & 0 \\ 0 & 0 & 0 & 0 & 0 \end{bmatrix}$$

$$\begin{bmatrix} 0 & 4 & 16 & 4 & 0 \\ \frac{16}{16} & 64 & 16 & 0 & 0 \\ 4 & 16 & 4 & 0 & 0 \\ 0 & 0 & 0 & 0 & 0 \\ 0 & 0 & 0 & 0 & 0 \end{bmatrix} \begin{bmatrix} 0 & 1 & 0 & 0 & 0 \\ 4 & 16 & 4 & 0 & 0 \\ 1 & 64 & \frac{7}{8} & 0 & 0 \\ 0 & 28 & 0 & 0 & 0 \\ 0 & 0 & 0 & 0 & 0 \end{bmatrix} \begin{bmatrix} 0 & 0 & 1 & 0 & 0 \\ 0 & 3 & 14 & 8 & 1 \\ 0 & 48 & 48 & 0 & 0 \\ 0 & 12 & 0 & 0 & 0 \\ 0 & 0 & 0 & 0 & 0 \end{bmatrix} \begin{bmatrix} 0 & 0 & 0 & 0 & 0 \\ 0 & \frac{23}{2} & \frac{17}{2} & 2 & 0 \\ 2 & \frac{28}{7} & \frac{64}{7} & \frac{4}{1} & 0 \\ 0 & 0 & 0 & 0 & 0 \\ 0 & 0 & 0 & 0 & 0 \end{bmatrix} \begin{bmatrix} 0 & 0 & 0 & 0 & 0 \\ 0 & 2 & \frac{10}{7} & 2 & 0 \\ 2 & 16 & 64 & 16 & 2 \\ 0 & 4 & 16 & 4 & 0 \\ 0 & 0 & 0 & 0 & 0 \end{bmatrix}$$

Table I. 5 × 5 stencils, scaled by 144 to convey relative size of entries while using the fewest digits, of the bi-4 coefficients of the construction for $\bar{T}$-nets. Each stencil lists the weight of the control nodes surrounding the $\bar{T}$-net.

14